\documentclass[10pt]{amsart}
\usepackage{amsfonts}
\usepackage{amssymb,latexsym}
\usepackage{amsmath}
\usepackage{amsthm}
\usepackage{amssymb}
\usepackage{enumerate}
\newtheorem{theorem}{Theorem}[section]
\newtheorem{lemma}[theorem]{Lemma}
\newtheorem{proposition}[theorem]{Proposition}

\theoremstyle{definition}

\theoremstyle{remark}
\newtheorem{remark}[theorem]{Remark}
\numberwithin{equation}{section}

\begin{document}
\setlength{\baselineskip}{1.2\baselineskip}
\title [A Minimal Value Problem]
{A Minimal Value Problem and the Prescribed $\sigma_2$ Curvature
Measure Problem}
\author{Chuanqiang Chen}
\address{Department of Mathematics\\
         University of Science and Technology of China\\
         Hefei, 230026, Anhui Province, CHINA}
\email{cqchen@mail.ustc.edu.cn}
\thanks{Research of the author was supported by Grant 10871187 from the National Natural Science
Foundation of China.}
\maketitle

\begin{abstract}
In this paper, we consider a minimal value problem and obtain an
algebraic inequality. As an application, we prove the $C^2$ a priori
estimate for a class of prescribed $\sigma_2$ curvature measure
equations, which generalizes the results of the $\sigma_2$ case  in
Guan-Li-Li\cite{GLL11} by a different method.
\end{abstract}

\section{Introduction}

In this paper, we consider the minimal value problem of the
following polynomial,
\begin{equation}\label{1.1}
f(x_1 , \cdots ,x_N) =  - b\sum\limits_{i = 1}^N {x_i } -
\sum\limits_{1 \le i < j \le N} {x_i x_j },
\end{equation}
with
\begin{equation}\label{1.2}
\sum\limits_{i = 1}^N {a_i x_i }  + C=0,
\end{equation}
 where $N\geq 2$ is an integer, $b$, $C$, $a_i$ are constants, and $\textbf{a}=(a_1,a_2, \cdots, a_N)$
satisfies the following two conditions
\begin{eqnarray}
\label{1.3}&&\sum\limits_{j = 1}^N {a_j } \neq 0, \\
\label{1.4}&&(\sum\limits_{j = 1}^N {a_j } )^2  - (N -
1)\sum\limits_{j = 1}^N {a_j ^2 } > 0.
\end{eqnarray}

\begin{lemma}\label{lem1.1}
Under the above assumptions, we have
\begin{equation}\label{1.5}
f(x_1 , \cdots ,x_N) \geq \frac{{b^2 [(\sum\limits_{j = 1}^N {a_j }
)^2  - N\sum\limits_{j = 1}^N {a_j ^2 } ] + 2 b C \sum\limits_{j =
1}^N {a_j } - (N - 1)C^2 }}{{2[(\sum\limits_{j = 1}^N {a_j } )^2  -
(N - 1)\sum\limits_{j = 1}^N {a_j ^2 } ]}}.
\end{equation}
\end{lemma}
Now, we recall the definition and some basic properties of
elementary symmetric functions. For any $k = 1, 2,\cdots, n,$ we set
$$
\sigma_k(\lambda) = \sum _{1 \le i_1 < i_2 <\cdots<i_k\leq
n}\lambda_{i_1}\lambda_{i_2}\cdots\lambda_{i_k},
 \qquad \text {for any} \quad\lambda=(\lambda_1,\cdots,\lambda_n)\in {\Bbb R}^n.
$$
We also set $\sigma_0=1$ and $\sigma_k =0$ for $k>n$.

In some sense, Lemma \ref{lem1.1} is corresponding to the convexity
of $\sigma_2$ operator (see Remark \ref{rmk4.2}).

As an application, we consider the prescribed curvature measures
problem, which is an important issue in convex geometry. There is a
vast literature devoted to the study of this type of problems(see
\cite{GM04, GLM09, GLL11} and references therein).

Let $X:M \rightarrow \mathbb{R}^{n+1}$ be a closed star-shaped
hypersurface in $\mathbb{R}^{n+1}$. The corresponding prescribed
$\sigma_k$ curvature measure equation is
\begin{equation}\label{1.6}
\sigma _k (\lambda(h_{ij}) ) = |X|^{-(n+1)} g(\frac{X}{|X|})
\left\langle {X,N} \right\rangle,
\end{equation}
where $N$ is the unit outer normal of $M$, $g \in C^2(\mathbb{S}^n)$
is a positive function and $h_{ij}$ is the second fundamental form
of $M$. Recall that the Garding's cone is defined as
\begin{equation}
\Gamma _k  = \{ \lambda  \in \mathbb{R}^n :\sigma _i (\lambda ) >
0,\forall 1 \le i \le k\}. \notag
\end{equation}
$M$ is called $k$-convex if its principal curvature vector
$\lambda(h_{ij}) \in \Gamma _k $, and the corresponding solution to
\eqref{1.6} is a $k$-admissible solution.

The existence theorem for $\sigma_k$ curvature measure problems ($1
\leq k \leq n$) has been solved. As the case $k = n$ is the
Alexandrov problem which was completely settled. Recently,
Guan-Lin-Ma \cite{GLM09} proved the existence theorem to the convex
solution case for $1 \leq k < n$ and the admissible solution for
$k=1$, and Guan-Li-Li \cite{GLL11} finished the admissible solution
case for $1 < k < n$ by establishing the key $C^2$ a priori
estimates.

The key $C^2$ a priori estimates in \cite{GLL11} is as follows.
\begin{theorem}\label{th1.2}
If $M$ satisfies equation \eqref{1.6} and $2 \leq k \leq n$, then
there exists a constant $C$ depending only on $n$, $k$,
$\min_{\mathbb{S}^n} g$, $||g||_{C^{1}}$ , and $||g||_{C^{2}}$, such
that
\begin{equation}\label{1.7}
\max_{M} \sigma_1 \leq C.
\end{equation}
\end{theorem}

And Guan-Li-Li \cite{GLL11} get the following existence theorem of
the admissible solution to \eqref{1.6}.
\begin{theorem}\label{th1.3}
Let $n\geq 2$ and $1 \leq k \leq n$. Suppose $g \in
C^2(\mathbb{S}^n)$ and $g > 0$. Then there exists a unique
$k$-convex star-shaped hypersurface $M \in C^{3,\alpha}$, $\forall
\alpha \in (0, 1)$ such that it satisfies \eqref{1.6}.
\end{theorem}
In \cite{GLL11}, Guan-Li-Li raised a question regarding global $C^2$
estimates for general curvature equations. That is, suppose $M
\subset \mathbb{R}^{n+1}$ is a compact smooth hypersurface
satisfying equation
\begin{equation}
\sigma _k (\kappa(X)) = g(X,N), \kappa(X) \in \Gamma_k, \forall X
\in M,
\end{equation}
where $g \in C^2(\mathbb{R}^{n+1} \times \mathbb{S}^n)$  is a
general positive function and $k\geq 2$. Suppose there is a priori
$C^1$ bound of $M$, can one conclude a $C^2$ a priori bound of $M$
in terms of $C^1$ norm of $M$, $g$, $n$, $k$?

In this paper, we give an affirmative answer to the following
equations,
\begin{equation}\label{1.9}
\sigma _2 (h_{ij} ) = \varphi(X) \left\langle {X,N}
\right\rangle^{\alpha}.
\end{equation}
where $\varphi\in C^2(\mathbb{R}^{n+1})$ and $\alpha \in
(-\infty,1+\delta)$ with $\delta=\inf_{M}{\frac{\left\langle {X,N}
\right\rangle^2}{|X|^2}}$. And we get the following $C^2$ a priori
estimate using Lemma \ref{lem1.1}.
\begin{theorem}\label{th1.4}
Let $n\geq 2$ and suppose $\varphi\in C^2(\mathbb{R}^{n+1})$ with
$\varphi > 0$. If $M$ is 2-convex and satisfies equation \eqref{1.9}
with the $C^1$ a priori estimate holding, then for $\alpha \in
(-\infty,1+\delta)$ with $\delta=\inf_{M}{\frac{\left\langle {X,N}
\right\rangle^2}{|X|^2}}$, there is a constant $C$ depending only on
$n$, $\alpha$, $||\varphi||_{C^{2}}$, $||\frac{1}{\varphi}||_{C^0}$
and the $C^1$ norm of $M$, such that
\begin{equation}\label{1.7}
\sup_{|\xi|=1}\sum\limits_{ij = 1}^n {h_{ij} \xi _i \xi _j} \leq C.
\end{equation}
\end{theorem}
When $\alpha =1$ and $\varphi(X)=|X|^{-(n+1)} g(\frac{X}{|X|})$, the
$C^1$ a priori estimate to \eqref{1.9} was proved in \cite{GLM09}.
So the $\alpha =1$ case of Theorem \ref{th1.4} is corresponding to
Theorem \ref{th1.2}.

The rest of this paper is organized as follows. In Section 2, we
prove Lemma \ref{lem1.1} by two methods. In Section 3, we prove
Theorem \ref{th1.4} using Lemma \ref{lem1.1}. And in Section 4, we
give some remarks.

\textbf{Acknowledgement} The author would like to express sincere
gratitude to Prof. Xi-Nan Ma for his encouragement and many
suggestions in this subject.

\section{Proof of Lemma \ref{lem1.1}}

In this section, we will prove lemma \ref{lem1.1} with two methods.
One is the eigenvector decomposition method, and the other is the
Lagrange method of multipliers. They are all very useful, and may be
used for other problems.

\textbf{Method 1. the eigenvector decomposition method.} We will
consider the eigenvalues and eigenvectors of the quadratic terms of
$f$.

By direct computations, we can know that the matrix
\begin{equation}\label{2.1}
\frac{1}{2}\left( {\begin{array}{*{20}c}
   0 & { - 1} &  \cdots  & { - 1}  \\
   { - 1} &  \ddots  &  \ddots  &  \vdots   \\
    \vdots  &  \ddots  &  \ddots  & { - 1}  \\
   { - 1} &  \cdots  & { - 1} & 0  \\
\end{array}} \right)
\end{equation}
has eigenvalues $-\frac{1}{2}(N-1)$ with a eigenvector $e_N  =
\frac{1}{{\sqrt N }}(1, \cdots ,1)$, and $\frac{1}{2}$ with
eigenvectors $e_1, \cdots, e_{N-1}$. In particular, we can assume
that $\{e_1, \cdots, e_N \}$ is an orthonormal basis of
$\mathbb{R}^N$ and
\begin{equation}\label{2.2}
\textbf{a} = (a_1 , \cdots ,a_N ) = \left\langle {\textbf{a},e_N }
\right\rangle e_N + d e_1  = \frac{{\sum\limits_{i = 1}^N {a_i }
}}{{\sqrt N }}e_N  + d e_1.
\end{equation}
So
\begin{equation}\label{2.3}
x = (x_1 , \cdots ,x_N ) = \sum\limits_{i = 1}^N {\left\langle
{x,e_i } \right\rangle e_i },
\end{equation}
and
\begin{equation}\label{2.4}
d^2  = \sum\limits_{i = 1}^N {a_i ^2 }  - \frac{{(\sum\limits_{i =
1}^N {a_i } )^2 }}{N}.
\end{equation}

When $d=0$, that is $a_1=a_2= \cdots =a_N \ne 0$, we can get
\begin{equation}
0 = \sum\limits_{i = 1}^N {a_i x_i }  + C = \left\langle
{x,\textbf{a}} \right\rangle  + C =\sqrt N a_1 \left\langle {x,e_N }
\right\rangle + C, \notag
\end{equation}
so
\begin{equation}\label{2.5}
\left\langle {x,e_N } \right\rangle = -\frac{C}{\sqrt N a_1}.
\end{equation}
Hence
\begin{eqnarray}\label{2.6}
f(x_1 , \cdots ,x_N ) =&&  - b\sum\limits_{i = 1}^N {x_i }  -
\sum\limits_{1 \le i < j \le N} {x_i x_j } \notag \\
= && - \sqrt N b\left\langle {x,e_N } \right\rangle  +
\frac{1}{2}\sum\limits_{i = 1}^{N - 1} {\left\langle {x,e_i }
\right\rangle ^2 }  - \frac{1}{2}(N - 1)\left\langle {x,e_N }
\right\rangle ^2 \notag \\
\geq && \frac{bC}{a_1} - \frac{1}{2}(N - 1)\frac{C^2}{N a_1^2}.
\end{eqnarray}
Hence \eqref{1.5} holds in this case.

When $d \ne 0$, we can get
\begin{eqnarray}\label{2.7}
0 = &&\sum\limits_{i = 1}^N {a_i x_i }  + C = \left\langle {x,a}
\right\rangle  + C \notag \\
=&& \frac{{\sum\limits_{i = 1}^N {a_i } }}{{\sqrt N }}\left\langle
{x,e_N } \right\rangle  + d\left\langle {x,e_1 } \right\rangle  + C,
\end{eqnarray}
so
\begin{equation}\label{2.8}
\left\langle {x,e_1 } \right\rangle  =  - \frac{{\sum\limits_{i =
1}^N {a_i } }}{{\sqrt N d}}\left\langle {x,e_N } \right\rangle  -
\frac{C}{d}.
\end{equation}
Then we can get
\begin{eqnarray}\label{2.9}
f(x_1 , \cdots ,x_N ) =&&  - b\sum\limits_{i = 1}^N {x_i }  -
\sum\limits_{1 \le i < j \le N} {x_i x_j } \notag \\
= && - \sqrt N b\left\langle {x,e_N } \right\rangle  +
\frac{1}{2}\sum\limits_{i = 1}^{N - 1} {\left\langle {x,e_i }
\right\rangle ^2 }  - \frac{1}{2}(N - 1)\left\langle {x,e_N }
\right\rangle ^2 \notag \\
=&&\frac{1}{2}\Bigg[\frac{{\sum\limits_{i = 1}^N {a_i } }}{{\sqrt N
d}}\left\langle {x,e_N } \right\rangle  + \frac{C}{d}\Bigg]^2 +
\frac{1}{2}\sum\limits_{i = 2}^{N - 1} {\left\langle {x,e_i }
\right\rangle ^2 }- \frac{1}{2}(N - 1)\left\langle {x,e_N }
\right\rangle ^2  - \sqrt N b\left\langle {x,e_N } \right\rangle
\notag \\
=&& \frac{1}{2}\Bigg[\frac{{(\sum\limits_{i = 1}^N {a_i } )^2  - (N
- 1)Nd^2 }}{{Nd^2 }}\left\langle {x,e_N } \right\rangle ^2  +
2\frac{{\sum\limits_{i = 1}^N {a_i } C - Nbd^2 }}{{\sqrt N d^2
}}\left\langle {x,e_N } \right\rangle  + \frac{{C^2
}}{{d^2 }}\Bigg]   \\
&& +\frac{1}{2}\sum\limits_{i = 2}^{N - 1} {\left\langle {x,e_i }
\right\rangle ^2 }.\notag
\end{eqnarray}
By \eqref{1.4},
\begin{eqnarray} \label{2.10}
(\sum\limits_{i = 1}^N {a_i } )^2  - (N - 1)Nd^2=N
\Bigg[(\sum\limits_{j = 1}^N {a_j } )^2  - (N - 1)\sum\limits_{j =
1}^N {a_j ^2 }\Bigg]>0,
\end{eqnarray}
so $f$ has a minimal value. And
\begin{eqnarray} \label{2.11}
f(x_1 , \cdots ,x_N ) =&& \frac{1}{2}\frac{{(\sum\limits_{i = 1}^N
{a_i } )^2  - (N - 1)Nd^2 }}{{Nd^2 }}\Bigg[\left\langle {x,e_N }
\right\rangle + \frac{{\sqrt N (\sum\limits_{i = 1}^N {a_i } C -
Nbd^2 )}}{{(\sum\limits_{i = 1}^N {a_i } )^2  - (N - 1)Nd^2
}}\Bigg]^2 \notag
\\
&&- \frac{{(\sum\limits_{i = 1}^N {a_i } C - bNd^2 )^2 }}{{2d^2
[(\sum\limits_{i = 1}^N {a_i } )^2  - (N - 1)Nd^2 ]}} + \frac{{C^2
}}{{2d^2 }} +\frac{1}{2}\sum\limits_{i = 2}^{N - 1}
{\left\langle {x,e_i } \right\rangle ^2 }  \notag \\
\ge&&  - \frac{{(\sum\limits_{i = 1}^N {a_i } C - bNd^2 )^2 }}{{2d^2
[(\sum\limits_{i = 1}^N {a_i } )^2  - (N - 1)Nd^2 ]}} + \frac{{C^2
}}{{2d^2 }} \notag \\
=&&\frac{{b^2 [(\sum\limits_{j = 1}^N {a_j } )^2  - N\sum\limits_{j
= 1}^N {a_j ^2 } ] + 2 b C \sum\limits_{j = 1}^N {a_j } - (N - 1)C^2
}}{{2[(\sum\limits_{j = 1}^N {a_j } )^2  - (N - 1)\sum\limits_{j =
1}^N {a_j ^2 } ]}}.
\end{eqnarray}
So Lemma \ref{lem1.1} holds.

\textbf{Method 2. the Lagrange method of multipliers.} Following the
same argument (only consider the quadratic terms of $\left\langle
{x,e_i } \right\rangle$), we can know $f$ has a minimal value by
\eqref{2.10}. Now we consider
\begin{equation}\label{2.12}
f(x_1 , \cdots ,x_N ,\mu ) =  - b\sum\limits_{i = 1}^N {x_i } -
\sum\limits_{1 \le i < j \le N} {x_i x_j }  + \mu [\sum\limits_{i =
1}^N {a_i x_i }  + C],
\end{equation}
so
\begin{eqnarray}
\label{2.13}&&0=f_{x_i} =  - b - \sum\limits_{j \ne i} { x_j } + \mu
a_i, \quad\forall \quad i=1, 2,
\cdots, N.\\
\label{2.14}&&0=f_{\mu} = \sum\limits_{i = 1}^N {a_i x_i }  + C.
\end{eqnarray}

By direct computations, we can get the minimal point,
\begin{eqnarray}
\label{2.15}&&\mu ^0  = \frac{{b\sum\limits_{j = 1}^N {a_j }  - (N -
1)C}}{{[\sum\limits_{j = 1}^N {a_j } ]^2  - (N - 1)\sum\limits_{j =
1}^N {a_j ^2 } }}, \\
\label{2.16}&&x_i^0  =  - \frac{1}{{N - 1}}b + \frac{{\mu ^0 }}{{N -
1}}\sum\limits_{j = 1}^N {a_j }  - \mu ^0 a_i, \quad\forall \quad
i=1, 2, \cdots, N.
\end{eqnarray}
So the minimum is
\begin{eqnarray}
f(x_1^0 , \cdots ,x_n^0 ,\mu ^0 ) =&&  - b\sum\limits_{i = 1}^N
{x_i^0 }  - \sum\limits_{1 \le i < j \le N} {x_i^0 x_j^0 } \notag \\
=&&  - b\sum\limits_{i = 1}^N {x_i^0 }  - \frac{1}{2}\sum\limits_{i
= 1}^N {x_i^0 \sum\limits_{j \ne i} {x_j^0 } }\notag \\
=&&  - b\sum\limits_{i = 1}^N {x_i^0 }  - \frac{1}{2}\sum\limits_{i
= 1}^N {x_i^0 [ - b + \mu ^0 a_i ]}\notag \\
=&&   - \frac{b}{2}\sum\limits_{i = 1}^N {x_i^0 }  - \mu ^0
\sum\limits_{i = 1}^N {a_i x_i^0 }\notag \\
=&& - \frac{b}{2}[ - \frac{N}{{N - 1}}b + \frac{{\mu ^0 }}{{N -
1}}\sum\limits_{j = 1}^N {a_j } ] + \mu ^0 C\notag \\
=&&\frac{{b^2 [(\sum\limits_{j = 1}^N {a_j } )^2  - N\sum\limits_{j
= 1}^N {a_j ^2 } ] + 2 b C \sum\limits_{j = 1}^N {a_j } - (N - 1)C^2
}}{{2[(\sum\limits_{j = 1}^N {a_j } )^2  - (N - 1)\sum\limits_{j =
1}^N {a_j ^2 } ]}}.
\end{eqnarray}

Now we finish the proof of Lemma \ref{lem1.1}.

\section{Proof of Theorem \ref{th1.4}}

First we collect some well-known properties of elementary symmetric
functions, which will be used in the proof of Theorem \ref{th1.4}.
\begin{proposition}\label{prop3.1}
Let $\lambda=(\lambda_1,\dots,\lambda_n)\in\mathbb{R}^n$ and $k =0,
1,...,n,$ then
\begin{align}
\label{3.1}&\sigma_k(\lambda)=\sigma_k(\lambda|i)+\lambda_i\sigma_{k-1}(\lambda|i), \quad\forall \,1\leq i\leq n,&\\
\label{3.2}&\sum_i^n \lambda_i\sigma_{k-1}(\lambda|i)=k\sigma_{k}(\lambda),&\\
\label{3.3}&\sum_i^n\sigma_{k}(\lambda|i)=(n-k)\sigma_{k}(\lambda).&
\end{align}
\end{proposition}
And we also have,
\begin{proposition} \label{prop3.2}
Suppose $W=(W_{ij})$ is diagonal, and $m$ ($1 \leq m \leq n$) is
positive integer, then
\begin{align}\label{3.4}
\frac{{\partial \sigma _m (W)}} {{\partial W_{ij} }} = \begin{cases}
\sigma _{m - 1} (W\left| i \right.), &\text{if } i = j, \\
0, &\text{if } i \ne j.
\end{cases}
\end{align}
and
\begin{align}\label{3.5}
\frac{{\partial ^2 \sigma _m (W)}} {{\partial W_{ij} \partial W_{kl}
}} =\begin{cases}
\sigma _{m - 2} (W\left| {ik} \right.), &\text{if } i = j,k = l,i \ne k,\\
- \sigma _{m - 2} (W\left| {ik} \right.), &\text{if } i = l,j = k,i \ne j,\\
0, &\text{otherwise }.
\end{cases}
\end{align}
\end{proposition}

Now, we recall some relevant geometric quantities of a smooth closed
hypersurface $M \subset \mathbb{R}^{n+1}$. Throughout the paper,
repeated indices denote summation and we assume the origin is inside
the body enclosed by $M$.

Let $M^n$ be an immersed hypersurface in $\mathbb{R}^{n+1}$. For $X
\in M \subset \mathbb{R}^{n+1}$, choose local normal coordinates in
$\mathbb{R}^{n+1}$, such that $\{ \frac{\partial}{\partial x_1},
\cdots,\frac{\partial}{\partial x_n}\}$ are tangent to $M$ and
$\partial_{n+1}$ is the unit outer normal of the hypersurface. We
sometimes denote $\partial_i=\frac{\partial}{\partial x_i}$ and also
use $N$ to denote the unit outer normal $\partial_{n+1}$. We use
lower indices to denote covariant derivatives with respect to the
induced metric.

For the immersion $X$, the second fundamental form is the symmetric
(2, 0)-tensor given by the matrix $\{h_{ij}\}$,
\begin{equation}\label{3.6}
h_{ij}=\left\langle {\partial_i X, \partial_j N} \right\rangle.
\end{equation}

Recall the following identities:
\begin{eqnarray}
\label{3.7}&&X_{ij} = -h_{ij}N  \quad \text{ (Gauss formula),} \\
\label{3.8}&&N_i = h_{ij}\partial_j \quad \text{(Weigarten equation), }\\
\label{3.9}&&h_{ijk} = h_{ikj} \quad \text{(Codazzi formula),} \\
\label{3.10}&&R_{ijkl} = h_{ik}h_{jl}- h_{il}h_{jk} \quad
\text{(Gauss equation)},
\end{eqnarray}
where $R_{ijkl}$ is the (4, 0)-Riemannian curvature tensor. We also
have
\begin{equation}\label{3.11}
h_{ijkl} =h_{klij} + (h_{mj}h_{il}- h_{ml}h_{ij})h_{mk} +
(h_{mj}h_{kl}-h_{ml}h_{kj} )h_{mi}.
\end{equation}

We will establish the curvature estimates to equation \eqref{1.9}.
Let
\begin{equation}\label{3.12}
\phi(X,\xi)  = \log (h_{kl}\xi_k \xi_l)(X)  - \log \left\langle
{X,N} \right\rangle,
\end{equation}
where $\xi \in \mathbb{S}^n$. Suppose that $\phi(X,\xi)$ attains its
maximum at some $X_0 \in M$ and $\xi_0 \in \mathbb{S}^n$. We may
assume $\xi_0 $ is $e_1$ and the other directions $e_2, \cdots, e_n$
can be chosen such that $\{e_1, e_2, \cdots , e_n\}$ is a local
orthonormal frame near $X_0$ and $h_{ij} (X_0)$ is diagonal. Then
the function
\begin{equation}\label{3.13}
\phi(X) = \log h_{11}  - \log \left\langle {X,N} \right\rangle.
\end{equation}
attains its maximum at $X_0 \in M$. All the following computations
is at the point $X_0$, and we denote $\lambda_i =h_{ii}$. So we have
\begin{equation}\label{3.14}
\phi _i  = \frac{{h_{11i} }}{{h_{11} }} - \frac{{\left\langle {X,N}
\right\rangle _i }}{{\left\langle {X,N} \right\rangle }} = 0,
\end{equation}
so
\begin{equation}\label{3.15}
\frac{{h_{11i} }}{{h_{11} }} = \frac{{\left\langle {X,N}
\right\rangle _i }}{{\left\langle {X,N} \right\rangle }}  =
\frac{{h_{ii}\left\langle {X,e_i} \right\rangle }}{{\left\langle
{X,N} \right\rangle }}.
\end{equation}
And
\begin{equation}\label{3.16}
\phi _{ii}  = \frac{{h_{11ii} }}{{h_{11} }} - \frac{{h_{11i} ^2
}}{{h_{11} ^2 }} - \frac{{\left\langle {X,N} \right\rangle _{ii}
}}{{\left\langle {X,N} \right\rangle }} + \frac{{\left\langle {X,N}
\right\rangle _i ^2 }}{{\left\langle {X,N} \right\rangle ^2 }} =
\frac{{h_{11ii} }}{{h_{11} }} - \frac{{\left\langle {X,N}
\right\rangle _{ii} }}{{\left\langle {X,N} \right\rangle }},
\end{equation}
hence
\begin{equation}\label{3.17}
0 \ge \sum\limits_{i = 1}^n {F^{ii} \phi _{ii} }  = \frac{1}{{h_{11}
}}\sum\limits_{i = 1}^n {F^{ii} h_{11ii} }  - \frac{1}{{\left\langle
{X,N} \right\rangle }}\sum\limits_{i = 1}^n {F^{ii} \left\langle
{X,N} \right\rangle _{ii} }.
\end{equation}
By \eqref{3.11}, we have
\begin{equation}\label{3.18}
h_{11ii}  = h_{ii11}  + h_{11} ^2 h_{ii}  - h_{11} h_{ii} ^2,
\end{equation}
so we obtain
\begin{eqnarray}\label{3.19}
\sum\limits_{i = 1}^n {F^{ii} h_{11ii} }  =&& \sum\limits_{i = 1}^n
{F^{ii} h_{ii11} }  + h_{11} ^2 \sum\limits_{i = 1}^n {F^{ii} h_{ii}
}  - h_{11} \sum\limits_{i = 1}^n {F^{ii} h_{ii} ^2 } \notag \\
=&& (\varphi \left\langle {X,N} \right\rangle^{\alpha} )_{11}  -
\sum\limits_{ijkl = 1}^n {F^{ij,kl} h_{ij1} h_{kl1} }  + 2 h_{11} ^2
\varphi \left\langle {X,N} \right\rangle^{\alpha}  - h_{11}
\sum\limits_{i = 1}^n {F^{ii} h_{ii} ^2 }\notag \\
=&& \varphi _{11} \left\langle {X,N} \right\rangle^{\alpha}  +
2{\alpha}\varphi _1 \left\langle {X,N} \right\rangle^{\alpha-1}
\left\langle {X,N} \right\rangle_1 + \alpha (\alpha - 1) \varphi
\left\langle {X,N} \right\rangle^{\alpha-2}\left\langle {X,N} \right\rangle_1^2 \notag \\
&& + \alpha \varphi \left\langle {X,N}
\right\rangle^{\alpha-1}[h_{11l} \left\langle {X,e_l } \right\rangle
+ h_{11} - h_{11} ^2 \left\langle {X,N}
\right\rangle ]\notag \\
&&- \sum\limits_{ijkl = 1}^n {F^{ij,kl} h_{ij1} h_{kl1} }  + 2h_{11}
^2 \varphi \left\langle {X,N} \right\rangle^{\alpha} - h_{11}
\sum\limits_{i = 1}^n {F^{ii} h_{ii} ^2 }.
\end{eqnarray}
And
\begin{eqnarray}\label{3.20}
\sum\limits_{i = 1}^n {F^{ii} \left\langle {X,N} \right\rangle _{ii}
}  =&& \sum\limits_{i = 1}^n {F^{ii} [h_{iil} \left\langle {X,e_l }
\right\rangle  + h_{ii}  - h_{ii} ^2 \left\langle {X,N}
\right\rangle ]} \notag \\
= &&(\varphi \left\langle {X,N} \right\rangle^{\alpha} )_l
\left\langle {X,e_l } \right\rangle  + 2\varphi \left\langle {X,N}
\right\rangle^{\alpha} - \left\langle {X,N} \right\rangle
\sum\limits_{i = 1}^n {F^{ii} h_{ii} ^2 }\notag \\
= &&\varphi_l \left\langle {X,N} \right\rangle^{\alpha} \left\langle
{X,e_l } \right\rangle +\alpha\varphi \left\langle {X,N}
\right\rangle^{\alpha-1}\left\langle {X,N} \right\rangle_l
\left\langle {X,e_l } \right\rangle \notag \\
&&+ 2\varphi \left\langle {X,N} \right\rangle^{\alpha} -
\left\langle {X,N} \right\rangle \sum\limits_{i = 1}^n {F^{ii}
h_{ii} ^2 }.
\end{eqnarray}
By \eqref{3.17}, \eqref{3.19} and \eqref{3.20},
\begin{eqnarray}\label{3.21}
0 \ge&& \left\langle {X,N} \right\rangle \sum\limits_{i = 1}^n
{F^{ii} h_{11ii} }  - h_{11} \sum\limits_{i = 1}^n {F^{ii}
\left\langle {X,N} \right\rangle _{ii} } \notag \\
=&& (2-\alpha)\varphi \left\langle {X,N} \right\rangle ^{1+\alpha}
h_{11} ^2 + h_{11}  [2\alpha \varphi _1 \left\langle {X,e_1 }
\right\rangle - \varphi _l \left\langle {X,e_l } \right\rangle
- (2-\alpha)\varphi ]\left\langle {X,N} \right\rangle^{\alpha} \notag \\
&&+ \varphi _{11} \left\langle {X,N} \right\rangle ^{1+\alpha} +
\alpha (\alpha - 1) \varphi \left\langle {X,N}
\right\rangle^{\alpha-1}\left\langle {X,e_1} \right\rangle^2h_{11}
^2 - \left\langle {X,N} \right\rangle \sum\limits_{ijkl = 1}^n
{F^{ij,kl} h_{ij1} h_{kl1} }.
\end{eqnarray}

In the following, we will deal with the last term in \eqref{3.21},
that is
\begin{eqnarray}\label{3.22}
- \sum\limits_{ijkl = 1}^n {F^{ij,kl} h_{ij1} h_{kl1} } =&&
\sum\limits_{i \ne j} {h_{ij1} ^2 }  - \sum\limits_{i \ne j}{h_{ii1} h_{jj1} } \notag \\
=&& \sum\limits_{i \ne j} {h_{ij1} ^2 }  - 2h_{111} \sum\limits_{i =
2}^n {h_{ii1} }  - 2\sum\limits_{ 2 \le i < j \le n} { h_{ii1}
h_{jj1} }.
\end{eqnarray}

Differentiating \eqref{1.9} once, we can get,
\begin{equation}\label{3.23}
\sum\limits_{i = 1}^n {\sigma _1 (\lambda |i)h_{ii1} }  = (\varphi
\left\langle {X,N} \right\rangle^{\alpha} )_1.
\end{equation}
Denote
\begin{equation}\label{3.24}
f(x_2 , \cdots ,x_n) =  - h_{111} \sum\limits_{i = 2}^n {x_i }  -
\sum\limits_{2 \le i < j \le n} {x_i x_j },
\end{equation}
and we consider the minimal value problem of $f(x_2 , \cdots ,x_n)$
under \eqref{3.23}, that is
\begin{equation}\label{3.25}
\sum\limits_{i = 2}^n {\sigma _1 (\lambda |i) x_i }+[\sigma _1
(\lambda |1) h_{111} - (\varphi \left\langle {X,N}
\right\rangle^{\alpha} )_1] = 0.
\end{equation}
By direct computations, we get
\begin{eqnarray}
\label{3.26}&& \sum\limits_{i = 2}^n {\sigma _1 (\lambda
|i)}=(n-2)\sigma_1(\lambda)+\lambda_1>0, \\
\label{3.27}&& [\sum\limits_{i = 2}^n {\sigma _1 (\lambda
|i)}]^2-(n-2)\sum\limits_{i = 2}^n {\sigma _1 (\lambda
|i)}^2=(n-1)\lambda_1^2+2(n-2)\sigma_2(\lambda)>0.
\end{eqnarray}
So from lemma \ref{lem1.1}, we can get $f(x_2 , \cdots ,x_n)$ has a
minimum
\begin{equation}\label{3.28}
\frac{{h_{111}^2 \Bigg[(\sum\limits_{i = 2}^n {\sigma _1 (\lambda
|i)})^2 - (n-1)\sum\limits_{i = 2}^n {\sigma _1 (\lambda |i)^2}
\Bigg] + 2 h_{111} C \sum\limits_{i = 2}^n {\sigma _1 (\lambda |i)}
- (n- 2)C^2 }}{{2[(n-1)\lambda_1^2+2(n-2)\sigma_2(\lambda)]}},
\end{equation}
where
\begin{eqnarray}\label{3.29}
C=&&\sigma _1 (\lambda |1) h_{111} - (\varphi \left\langle {X,N}
\right\rangle^{\alpha} )_1   \notag \\
=&& \sigma _1 (\lambda |1) h_{111} - \varphi _1 \left\langle {X,N}
\right\rangle^{\alpha} -\alpha\varphi \left\langle
{X,N}\right\rangle^{\alpha-1} \left\langle {X,N}\right\rangle_1 \notag \\
=&& [\sigma _1 (\lambda |1)\lambda _1  -\alpha\sigma _2 (\lambda
)]\frac{{h_{111}}}{{h_{11}}}-\varphi_1 \left\langle
{X,N}\right\rangle^{\alpha}.
\end{eqnarray}
So
\begin{eqnarray}\label{3.30}
&&h_{111}^2 \Bigg[(\sum\limits_{i = 2}^n {\sigma _1 (\lambda |i)})^2
- (n-1)\sum\limits_{i = 2}^n {\sigma _1 (\lambda |i)^2} \Bigg] + 2
h_{111} C \sum\limits_{i = 2}^n {\sigma _1 (\lambda |i)}
- (n- 2)C^2 \notag \\
=&&(\frac{h_{111}}{h_{11}})^2 \Bigg[(\sum\limits_{i = 2}^n {\sigma
_1 (\lambda |i)})^2 - (n-1)\sum\limits_{i = 2}^n {\sigma _1 (\lambda
|i)^2}\Bigg] \lambda _1^2 \notag \\
&& + 2\frac{h_{111}}{h_{11}} \Bigg[ [\sigma _1 (\lambda |1)\lambda
_1  -\alpha\sigma _2 (\lambda )]\frac{{h_{111}}}{{h_{11}}}-\varphi_1
\left\langle {X,N}\right\rangle^{\alpha} \Bigg]\lambda
_1\sum\limits_{i =
2}^n {\sigma_1 (\lambda |i)}\notag \\
&&-(n-2)\Bigg[ [\sigma _1 (\lambda |1)\lambda _1  -\alpha\sigma _2
(\lambda )]\frac{{h_{111}}}{{h_{11}}}-\varphi_1 \left\langle
{X,N}\right\rangle^{\alpha} \Bigg]^2.
\end{eqnarray}
By direct computations, the coefficient of
$(\frac{h_{111}}{h_{11}})^2$ in \eqref{3.30} is
\begin{eqnarray}
&&{\lambda _1 }^2[(\sum\limits_{i = 2}^n {\sigma _1 (\lambda |i)})^2
-(n-1)\sum\limits_{i = 2}^n {\sigma _1 (\lambda |i)^2}] + 2[\sigma
_1 (\lambda |1)\lambda _1 -\alpha\sigma _2 (\lambda )]\lambda
_1\sum\limits_{i = 2}^n {\sigma_1 (\lambda |i)}\notag \\
&&-(n-2)[\sigma _1 (\lambda |1)\lambda _1 -\alpha\sigma _2 (\lambda
)]^2\notag \\
=&&{\lambda _1 }^2[(\sum\limits_{i = 2}^n {\sigma _1 (\lambda
|i)})^2 -(n-1)\sum\limits_{i = 2}^n {\sigma _1 (\lambda |i)^2} +
2\sigma _1 (\lambda |1)\lambda _1\sum\limits_{i = 2}^n {\sigma_1
(\lambda |i)}-(n-2)\sigma _1^2 (\lambda |1)]\notag \\
&&-2\alpha\sigma _2 (\lambda )\lambda_1[\sum\limits_{i = 2}^n {\sigma_1 (\lambda |i)}-(n-2)\sigma _1 (\lambda |1)]-(n-2)\alpha^2\sigma _2^2 (\lambda)\notag \\
=&&{\lambda _1 }^2[(\sum\limits_{i = 1}^n {\sigma _1 (\lambda
|i)})^2 -(n-1)\sum\limits_{i = 1}^n {\sigma _1 (\lambda |i)^2} ]\notag \\
&& -2(n-1)\alpha\sigma _2 (\lambda )\lambda_1^2-(n-2)\alpha^2\sigma _2^2 (\lambda)\notag \\
=&&{\lambda _1 }^2\Bigg[[(n-1) \sigma _1 (\lambda)]^2 -(n-1)[(n-1)
\sigma _1 (\lambda)^2 -2\sigma _2 (\lambda )] \Bigg]\notag \\
&& -2(n-1)\alpha\sigma _2 (\lambda )\lambda_1^2 -(n-2)\alpha^2\sigma
_2^2(\lambda)\notag \\
=&&2(1-\alpha)(n-1)\sigma _2 (\lambda )\lambda_1^2
-(n-2)\alpha^2\sigma _2^2(\lambda).\notag
\end{eqnarray}
So
\begin{eqnarray}
&&h_{111}^2 \Bigg[(\sum\limits_{i = 2}^n {\sigma _1 (\lambda |i)})^2
- (n-1)\sum\limits_{i = 2}^n {\sigma _1 (\lambda |i)^2} \Bigg] + 2
h_{111} C \sum\limits_{i = 2}^n {\sigma _1 (\lambda |i)}
- (n- 2)C^2 \notag \\
=&&[2(1-\alpha)(n-1)\sigma _2 (\lambda )\lambda_1^2
-(n-2)\alpha^2\sigma _2^2(\lambda)](\frac{h_{111}}{h_{11}})^2\notag \\
&& - 2\frac{h_{111}}{h_{11}} [(n-1)\lambda
_1^2+(n-2)\alpha\sigma_2(\lambda)]\varphi _1 \left\langle {X,N}
\right\rangle^{\alpha}-(n-2)\varphi _1^2 \left\langle {X,N} \right\rangle^{2\alpha}\notag \\
=&&2(1-\alpha)(n-1)\varphi\left\langle {X,N}
\right\rangle^{\alpha-2} \left\langle {X,N} \right\rangle_1^2
\lambda_1 ^2 - (n - 2)\alpha^2\varphi^2 \left\langle {X,e_1 }
\right\rangle ^2 \lambda _1 ^2 \notag \\
&& + 2\varphi _1 \left\langle {X,e_1 } \right\rangle [(n - 1)\lambda
_1 ^3  + (n - 2)\lambda _1 \alpha\varphi \left\langle {X,N}
\right\rangle^{\alpha} ] -(n-2)\varphi _1^2 \left\langle {X,N}
\right\rangle^{2\alpha}. \notag
\end{eqnarray}
Hence
\begin{eqnarray}
&&- \left\langle {X,N} \right\rangle \sum\limits_{ijkl = 1}^n
{F^{ij,kl}h_{ij1} h_{kl1} } \geq \left\langle {X,N} \right\rangle[-
2h_{111} \sum\limits_{i = 2}^n {h_{ii1} } -2\sum\limits_{ 2 \le i < j \le n} { h_{ii1} h_{jj1} }] \notag \\
\geq&& \left\langle {X,N} \right\rangle\frac{{h_{111}^2
\Bigg[(\sum\limits_{i = 2}^n {\sigma _1 (\lambda |i)})^2 -
(n-1)\sum\limits_{i = 2}^n {\sigma _1 (\lambda |i)^2} \Bigg] + 2
h_{111} C \sum\limits_{i = 2}^n {\sigma _1
(\lambda |i)} - (n- 2)C^2}}{{(n-1)\lambda_1^2+2(n-2)\sigma_2(\lambda)}}\notag \\
\geq&& (1 - \alpha) \varphi \left\langle {X,N}
\right\rangle^{\alpha-1}\left\langle {X,N}
\right\rangle_1^2[\frac{{2(n-1)\lambda_1 ^2 }} {{(n-1)\lambda_1^2 +
2(n-2)\sigma_2(\lambda) }}]-C_1h_{11}-C_2 \notag\\
=&& (1 - \alpha) \varphi \left\langle {X,N}
\right\rangle^{\alpha-1}\left\langle {X,e_1}
\right\rangle^2h_{11}^2[2-\frac{{4(n-2)\sigma_2(\lambda)}}
{{(n-1)\lambda_1^2 + 2(n-2)\sigma_2(\lambda) }}]-C_1h_{11}-C_2 \notag \\
\geq&& 2(1 - \alpha) \varphi \left\langle {X,N}
\right\rangle^{\alpha-1}\left\langle {X,e_1}
\right\rangle^2h_{11}^2-C_1h_{11}-C_3. \notag
\end{eqnarray}
From \eqref{3.21}, we can get
\begin{eqnarray}\label{3.31}
0 \ge&& (2-\alpha)\varphi \left\langle {X,N} \right\rangle
^{1+\alpha} h_{11} ^2 -Ch_{11}-C \notag \\
&& + \alpha (\alpha - 1) \varphi \left\langle {X,N}
\right\rangle^{\alpha-1}\left\langle {X,e_1} \right\rangle^2h_{11}
^2 - \left\langle {X,N} \right\rangle \sum\limits_{ijkl = 1}^n
{F^{ij,kl} h_{ij1} h_{kl1} } \notag \\
\ge&& (2-\alpha)\varphi\left\langle {X,N} \right\rangle ^{\alpha-1}
h_{11} ^2[\left\langle {X,N} \right\rangle ^2+(1-\alpha)\left\langle
{X,e_1} \right\rangle^2] -Ch_{11}-C.
\end{eqnarray}

When $\alpha \in (-\infty,1]$, we can get
\begin{equation}\label{3.32}
0 \ge (2-\alpha)\varphi\left\langle {X,N} \right\rangle ^{\alpha+1}
h_{11} ^2 -Ch_{11}-C,
\end{equation}
so Theorem \ref{th1.4} holds.

When $\alpha \in (1,1+\delta)$ with $\delta = \inf_{M}
{\frac{\left\langle {X,N} \right\rangle^2}{|X|^2}}$, we can get
\begin{eqnarray}
\left\langle {X,N} \right\rangle ^2+(1-\alpha)\left\langle {X,e_1}
\right\rangle^2 >&& \left\langle {X,N} \right\rangle
^2+(1-\alpha)|X|^2 \notag \\
=&&\left\langle {X,N} \right\rangle
^2-\delta|X|^2+(1+\delta-\alpha)|X|^2\notag \\
\geq&& (1+\delta-\alpha)|X|^2,
\end{eqnarray}
so we obtain
\begin{equation}\label{3.32}
0 \ge (2-\alpha)(1+\delta-\alpha)\varphi\left\langle {X,N}
\right\rangle ^{\alpha-1}|X|^2 h_{11} ^2 -Ch_{11}-C.
\end{equation}
So Theorem \ref{th1.4} holds.

\section{Some Remarks}

\begin{remark}
For the minimal value problem \eqref{1.1}, \eqref{1.2}, there is a
unique minimal point such that ''='' holds in \eqref{1.5}. So
\eqref{1.5} is optimal. And we can get the unique minimal point from
the proof.
\end{remark}

\begin{remark} \label{rmk4.2}
The minimal value problem \eqref{1.1}, \eqref{1.2} (with
$a_i=\sigma_1(\lambda|i)$, $i=2, \cdots, n$) is corresponding to the
convexity of $\sigma_2$. To be precise, the corresponding minimal
value problems is
\begin{equation}
f^2(x_2, \cdots ,x_n) =  - b\sum\limits_{i = 2}^n { x_i } -
\sum\limits_{2 \le i < j \le n} {x_i x_j },
\end{equation}
with
\begin{equation}
\sum\limits_{i = 2}^n {\sigma_{1}(\lambda|i) x_i }
+[\sigma_{1}(\lambda|1)b + C]=0.
\end{equation}
And we can get
\begin{equation}
f^2(x_2, \cdots ,x_n)\geq \frac{{
2(n-1)\sigma_2(\lambda)b^2+2(n-1)\lambda_1bC- (n- 2)C^2
}}{{2[(n-1)\lambda_1^2+2(n-2)\sigma_2(\lambda)]}}.
\end{equation}
 For general $\sigma_k$ ($k>2$), we guess
there are similar minimal value problems, that is
\begin{equation}
f^k(x_2, \cdots ,x_n) =  - b\sum\limits_{i = 2}^n
{\sigma_{k-2}(\lambda|1i)x_i } - \sum\limits_{2 \le i < j \le n}
{\sigma_{k-2}(\lambda|ij)x_i x_j },
\end{equation}
with
\begin{equation}
\sum\limits_{i = 2}^n {\sigma_{k-1}(\lambda|i) x_i }
+[\sigma_{k-1}(\lambda|1)b + C]=0.
\end{equation}
But we cannot get the minimum by the methods in Section 2. There are
some techniques which cannot be used here.
\end{remark}

\begin{remark} \label{rmk4.3}
For a general prescribed $\sigma_2$ curvature measure equation
\begin{equation}\label{4.3}
\sigma _2 (h_{ij} ) = \varphi(X) \left\langle {X,N}
\right\rangle^{\alpha},
\end{equation}
where $\alpha \in (-\infty,1+\delta)$, the $C^1$ a priori estimate
was proved in \cite{GLM09} when $\alpha=1$. In fact, the proof of
$C^1$ a priori estimate in \cite{GLM09} holds for $\alpha \in (0,1]$
with a small modification. For general $\alpha$, I don't know the
corresponding results because of my limitation of knowledge.
\end{remark}

\end{document}